\documentclass[a4paper,10.5pt]{amsart}

\usepackage{amsmath}
\usepackage{amssymb}
\usepackage{amsthm}
\usepackage{amscd}
\usepackage[dvips]{graphicx}
\usepackage{enumerate}
\usepackage{breqn}
\usepackage{cleveref}
\usepackage{enumerate}

\newtheorem{theorem}{Theorem}[section]

\numberwithin{equation}{section}
\newcommand{\real}{\mathbb{R}}

\def\eqn {\begin{equation}}
\def\eeqn {\end{equation}}
\def\C{{\mathbb C}}
\def\real{{\mathbb R}}

\def\pa{\partial}

\def\C{\mathcal C}

\def\F{\mathcal F}
\def\K{\mathcal K}

\def\S{\mathcal S}

\begin{document}
\title{Global Continuation and the Theory of Rotating Stars}
\author{Yilun Wu}
\address{Department of Mathematics, University of Oklahoma, Norman, OK 73019}
\email{allenwu@ou.edu}
\date{}

\begin{abstract}
This paper gives a condensed review of the history of solutions to the Euler-Poisson equations modeling equilibrium states of rotating stars and galaxies, leading to a recent result of Walter Strauss and the author. This result constructs a connected set of rotating star solutions for larger and larger rotation speed, so that the supports of the stars become unbounded if we assume an equation of state $p = \rho^\gamma$, $4/3<\gamma<2$. On the other hand, if $6/5<\gamma<4/3$, we show that either the supports of the stars become unbounded, or the density somewhere within the stars becomes unbounded. This is the first global continuation result for rotating stars that displays singularity formation within the solution set.
\end{abstract}

\maketitle

\section{A brief history on equilibrium of rotating fluids}
The equilibrium shape and density distribution of rotating fluids under self gravitation is a classical problem in mathematical physics with a long history. Such a fluid can be modeled by the Euler-Poisson equations, a system coupling perfect fluid with Newtonian gravity:
\begin{equation}\label{eq: full Euler-Poisson}
\begin{cases}
\rho _t + \nabla\cdot (\rho v)=0, \\
(\rho v)_t + \nabla\cdot(\rho v\otimes v) + \nabla p = \rho \nabla U, \\
U(x,t) = \int_{\real^3}\frac{\rho(x',t)}{|x-x'|}~dx'.
\end{cases}
\end{equation}
This system can be reduced to 
\begin{equation}
\begin{cases}\label{eq: EP1}
-\rho\omega^2(r)re_r +\nabla p =\rho\nabla U,\\
U(x,t) = \int_{\real^3}\frac{\rho(x',t)}{|x-x'|}~dx',
\end{cases}
\end{equation}
if we make the following assumptions
\begin{enumerate}
\item All functions are time independent.
\item $v = \omega(r) r e_\theta$.
\item $\rho$ is constant on the fluid domain (incompressible) or is a function of $r$ and $x_3$ only (compressible).
\end{enumerate}
In the above we have used the cylindrical coordinate $r=\sqrt{x_1^2+x_2^2}$, and the unit vectors $e_r = \frac1r(x_1,x_2,0)$ and $e_\theta = \frac1r(-x_2,x_1,0)$.
We require \eqref{eq: EP1} to hold on the fluid domain $\{\rho>0\}$. We also require the vacuum boundary condition:
\eqn
p = 0 \text{ on }\partial \{\rho > 0\}.
\eeqn
Newton essentially started thinking about near spherical solutions to \eqref{eq: EP1} soon after he discovered his law of gravity. Most of the early attempts in solving this problem involve trying the ansatz $\rho = \chi_D$, the characteristic function of a suitable smooth domain $D$ (thus describing an incompressible fluid), while setting $\omega(r)\equiv \omega_0$, a uniform rotation profile. Under these assumptions, one has:
\eqn\label{eq: EP constant density}
\begin{cases}
\nabla \left(-\frac12\omega_0^2 (x_1^2+x_2^2) - U+p\right)=0 \quad &\text{ on }D,\\
U = \frac1{|x|}*\chi_D \quad &\text{ on }\real^3,\\
p=0 \quad &\text{ on }\partial D.
\end{cases}
\eeqn
If we assume $D$ has only one connected component, \eqref{eq: EP constant density} essentially requires
\eqn\label{eq: domain equation}
\frac12\omega_0^2 (x_1^2+x_2^2) + \frac1{|x|}*\chi_D = \text{constant} \text{ on }\partial D.
\eeqn
Strictly speaking, one should also check that $p\ge 0$ on $D$, but this can be easily verified at the end and is omitted from the following discussion. One therefore just needs to find a domain $D$ for which \eqref{eq: domain equation} holds.
The first well-known exact solution of this sort is due to MacLaurin in the eighteenth century. He uses the formula for the Newtonian potential of an ellipsoid (see, for example, section VII.6 in \cite{kellogg1953foundations}): if $D$ is an ellipsoid $\left\{\frac{x_1^2}{a^2}+\frac{x_2^2}{b^2}+\frac{x_3^2}{c^2}\le 1\right\}$, then 
\eqn
\frac{1}{|x|}*\chi_D = L_0(a,b,c)- L_1(a,b,c)x_1^2-L_2(a,b,c)x_2^2-L_3(a,b,c)x_3^2
\eeqn
for $x\in \overline{D}$, where
\eqn
L_0(a,b,c)=\pi abc \int_0^\infty\frac{ds}{\sqrt{(a^2+s)(b^2+s)(c^2+s)}},
\eeqn
\eqn
L_1(a,b,c)=\pi abc\int_0^\infty \frac{ds}{(a^2+s)\sqrt{(a^2+s)(b^2+s)(c^2+s)}},
\eeqn
\eqn
L_2(a,b,c)=\pi abc\int_0^\infty \frac{ds}{(b^2+s)\sqrt{(a^2+s)(b^2+s)(c^2+s)}},
\eeqn
\eqn
L_3(a,b,c)=\pi abc\int_0^\infty \frac{ds}{(c^2+s)\sqrt{(a^2+s)(b^2+s)(c^2+s)}}.
\eeqn
Maclaurin looks for an axisymmetric ellipsoid for which $a=b$. Thus 
\eqn
\frac{1}{|x|}*\chi_D = L_0(a,a,c)-L_1(a,a,c)(x_1^2+x_2^2)-L_3(a,a,c)x_3^2
\eeqn
on $\overline D$. On $\partial D$, $x_3^2 = c^2 \left(1-\frac{x_1^2+x_2^2}{a^2}\right)$, thus 
\eqn
\frac{1}{|x|}*\chi_D = (L_0-c^2L_3)- \left(L_1-\frac{c^2}{a^2}L_3\right)(x_1^2+x_2^2).\eeqn
\eqref{eq: domain equation} now becomes 
\eqn
\left[\frac12\omega_0^2 -\left(L_1-\frac{c^2}{a^2}L_3\right)\right](x_1^2+x_2^2) = \text{constant on }\partial D,
\eeqn
which is satisfied if we take 
\eqn\label{eq: Mac ellip}
\frac12\omega_0^2  = L_1(a,a,c)-\frac{c^2}{a^2}L_3(a,a,c).
\eeqn
Let us consider solutions with fixed total mass. As the volume of the ellipsoid is $\pi a^2c$, let's set $a^2 c =1$ for simplicity. As a consequence, $a^3=\frac a c$. If we define $e=\frac a c$ to be the ellipticity of the ellipsoid, one can easily find
\eqn\label{eq: Mac int}
L_1(a,a,c) - \frac {c^2}{a^2}L_3(a,a,c) = \pi\int_0^\infty \frac{1}{(1+s)\sqrt{1+e^2s}}\left(\frac{1}{1+s}-\frac1{1+e^2s}\right)~ds.
\eeqn
By \eqref{eq: Mac ellip}, we get a Maclaurin ellipsoidal solution whenever the right hand side of \eqref{eq: Mac int} is nonnegative. This happens if and only if $e\ge 1$. In other words, the solutions are ``oblate". By the relation of $e$, $a$ and $c$ given above, when $e$ tends to infinity, $a$ tends to infinity and $c$ tends to zero. Thus we get a continuous set of solutions, so that the support of the fluid domain blows up along this set. It is interesting to note that the angular velocity $\omega_0$ does not blow up in this set, as the right hand side of \eqref{eq: Mac int} tends to zero as $e$ tends to infinity.

The Maclaurin ellipsoids present a simple example of a solution set that shows blow up behavior. Much of the recent progress made by Walter Strauss and the author is about constructing a similar solution set for the {\it compressible} Euler-Poisson equation, as I will show in the following. Nevertheless, the transition from the incompressible model to the compressible one happened rather slowly in history. In retrospect, this could be due to the fact that compressible solutions are much more difficult to construct and will need some serious input from modern PDE theory. 

To continue our discussion of the classical history, several other main events include: the discovery of other non-axisymmetric ellipsoidal solutions by Jacobi in the nineteenth century; the study of linear perturbations of these constant density ellipsoids by Poincar\'e, and the study of nonlinear perturbations by Lyapunov, both in the early twentieth centry. The solutions found by Poincar\'e and Lyapunov are both incompressible, and the density function $\rho$ is close to a constant on the fluid domain. A very nice account of the classical history of this problem, including discussions of the above mentioned works, can be found in Jardetzky \cite{jardetzky2013theories}.

To provide more realistic models of gaseous stars, people gradually turned to compressible gas dynamics, in which an equation of state $p=p(\rho)$ is prescribed to relate pressure directly to fluid density. In the late nineteenth century, Lane and Emden studied {\it non-rotating} star solutions under a power law: $p = C\rho^\gamma$ for some positive constants $C,\gamma$. Their work made a big impact on the study of stellar structure in astrophysics. Chandrasekhar \cite{chandrasekhar1939introduction} is a classical reference for this work. To explain in our current language, we take $\omega(r) = 0$ in \eqref{eq: EP1}, use the equation of state $p=C\rho^\gamma$, and divide the first equation by $\rho$ on the fluid domain. It follows that
\eqn
\nabla \left(\rho^{\gamma-1}-\frac{1}{|x|}*\rho\right)=0.
\eeqn
For simplicity of presentation, I have chosen a suitable $C$ so that the coefficient in front of $\rho^{\gamma-1}$ is $1$. By assuming the fluid domain has one connected component, we get 
\eqn\label{eq: LE1}
\rho^{\gamma-1}-\frac{1}{|x|}*\rho = \text{constant on }\{\rho>0\}.
\eeqn
Assume $\rho$ is radially symmetric and supported on a ball, \eqref{eq: LE1} is equivalent to 
\eqn
\Delta\left(\rho^{\gamma-1}-\frac{1}{|x|}*\rho\right)=0 \text{ on }\{\rho>0\}.
\eeqn
The reason is that any radially symmetric harmonic function on a ball centered at the origin is constant. Now letting $u = \rho^{\gamma-1}$, and using $\Delta^{-1} = -\frac{1}{4\pi|x|}*\cdot$ in $\real^3$, we get
\eqn\label{eq: LE2}
\Delta u + 4\pi u^{1/(\gamma-1)}=0 
\eeqn
on $\{u>0\}$. This is the well-known Lane-Emden equation. It is actually an ODE for radial solutions and can be treated as one. Alternatively, one can construct solutions using PDE methods, which provide a more uniform treatment even when the equation of state is not exactly a power law. Let us summarize the result by the following 
\begin{theorem}\label{thm: LE}
Consider \eqref{eq: LE2} on $\real^3$. The existence of compactly supported positive radial solution of \eqref{eq: LE2} depends on the exponent $\frac1{\gamma-1}$:
\begin{itemize}
\item If $0<\frac{1}{\gamma-1}<5$, $\frac{1}{\gamma-1}\ne 1$, then on any finite ball $B$ centered at the origin, there exists a unique positive radial solution to \eqref{eq: LE2}, such that it is continuous on $\overline B$ and $u=0$ on $\partial B$. 
\item If $\frac{1}{\gamma-1}=1$, then a positive solution with zero boundary value can only exist on the ball of radius $\frac{\sqrt\pi}{2}$, and the solution has the explicit formula $u(x) = C\frac{\sin(2\sqrt{\pi}|x|)}{|x|}$ for some positive constant $C$.
\item If $\frac{1}{\gamma-1}\ge5$, there is no positive solution with zero boundary value on any finite ball.
\end{itemize}
\end{theorem}
Here, existence of solutions for $0<\frac{1}{\gamma-1}<5$, $\frac{1}{\gamma-1}\ne 1$ is a special case of results in \cite{ambrosetti1973dual} and \cite{de1982priori}. A uniqueness proof can be found in \cite{gidas1979symmetry} or more generally in \cite{ni1985uniqueness}. The $\frac1{\gamma-1}=1, 5$ cases can be solved explicitly. The nonexistence result for $\frac1{\gamma-1}>5$ follows from the classical Pohozaev identity (See for example, section 9.4.2 in \cite{evans2010partial}). By Theorem \ref{thm: LE}, the range of $\gamma$ for existence is $\gamma>\frac65$.

After Lane and Emden's discovery, attempts were made to compute linear perturbations of these solutions to produce rotating stars (see \cite{chandrasekhar1967ellipsoidal}, for example), but the first nonlinear construction of exact rotating star solutions is due to Lichtenstein \cite{lichtenstein1933untersuchungen}. His result in our current language can be summarized as follows. As before, we divide the first equation in \eqref{eq: EP1} by $\rho$ on the fluid domain, and use the equation of state $p=C\rho^\gamma$. The equation can now be written as
\eqn
\nabla \left(\rho^{\gamma-1}-\frac1{|x|}*\rho - \int_0^r \omega^2(s)s~ds \right)=0 \text{ on }\{\rho>0\},
\eeqn
or 
\eqn\label{eq: rot star}
\rho^{\gamma-1}-\frac1{|x|}*\rho - \int_0^r \omega^2(s)s~ds = \text{constant on }\{\rho>0\} .
\eeqn
Lichtenstein's result can be described as
\begin{theorem}\label{thm: Lit}
Let $\frac65<\gamma<2$, and $\rho_0$ be a Lane-Emden solution. Let $\omega(r) = \kappa \omega_0(r)$, where $\omega_0(r)$ is any given smooth function. Then for each sufficiently small $\kappa$, there is a nonnegative compactly supported continuous function $\rho = \rho(\kappa)$ solving \eqref{eq: rot star}. The mapping $\kappa\mapsto \rho(\kappa)$ is continuous into a suitable function space, and $\rho(0)=\rho_0$.
\end{theorem}
Put more informally, he constructed a continuous curve of slowly rotating stars that are small perturbations of a given Lane-Emden solution. The range of $\gamma$ in this result is much more limited compared with the full range of the Lane-Emden solutions ($\gamma>\frac65$), but actually covers most types of gases relevant to astrophysics. It is worth noting that Lichtenstein's contruction is done in such a way that it is unclear whether the solutions obtained will have the same total mass as $\rho_0$. This is a topic that would be taken on by Walter and me later and would turn out to be an important issue for studying large deviations from $\rho_0$. 

Lichtenstein's work, unfortunately, did not make a significant impact on the rotating star literature. In retrospect, the reason for such limited impact may be twofold. To the astrophysics community, the construction of an exact solution may appear as a technical piece of mathematical curiosity, and would be less interesting than an actual calculation of the linear perturbation. On the mathematical side, Lichtenstein's proof of the result is not completely transparent with all the delicate estimates he needs to show convergence of his perturbation series. In fact, many years later, Heilig \cite{heilig1994lichtenstein} served to crystalize Lichtenstein's argument as an application of the implicit function theorem on a suitable function space. Even after Heilig's rework, Lichtenstein's result appears to remain relatively unknown to the mathematical community.

The next major event in the history is Auchmuty and Beals' work \cite{auchmuty1971variational}, which is the first result on rotating stars that does not require the rotation to be small. Their result can be described as follows:
\begin{theorem}
Let $\gamma>\frac43$, and $M>0$ be given, and let $\omega(r)$ be any given smooth function with sufficient decay at infinity. Then there exists a nonnegative compactly supported continuous function $\rho$ solving \eqref{eq: rot star}, such that $\int_{\real^3}\rho(x)~dx = M$.
\end{theorem}
This result has several advantages compared to Lichtenstein's. It covers a wide range of $\gamma$; it has a built in mass constraint; it does not require smallness of rotation. On the other hand, it has the disadvantage of requiring $\omega(r)$ to have a certain kind of decay at infinity. This drawback was partially removed by Li \cite{li1991uniformly}, who showed the same result for constant rotation profile $\omega(r) \equiv \omega_0 $. The method of \cite{auchmuty1971variational} is calculus of variations (energy minimization), and is completely different from Lichtenstein's perturbation method. \cite{auchmuty1971variational} made a big impact on the mathematical literature of rotating stars. Developing the variational techniques used in \cite{auchmuty1971variational}, Friedman and Turkington \cite{friedman1981existence}, Li \cite{li1991uniformly}, McCann \cite{mccann2006stable}, Wu \cite{wu2015rotating} and Wu \cite{wu2016existence}  proved existence results in various more general setups. Caffarelli and Friedman \cite{caffarelli1980shape}, Friedman and Turkington \cite{friedman1980asymptotic}, Chanillo and Li \cite{chanillo1994diameters} studied qualitative properties and bounds on the size of the support of the variational solutions. Luo and Smoller \cite{luo2009existence} proved a conditional nonlinear stability result using the variational method.

\section{Revival of the perturbation method}
By the time I went to Brown University as a postdoc working with Walter, my knowledge of the rotating star literature is pretty much dominated by the variational approach. We were not even aware of Lichtenstein's work which had been published some eighty years ago. At that point, Walter raised the interesting question of studying the continuity of the set of rotating star solutions, and whether certain forms of blow up may appear as one globally continues along the solution set. This is a natural analog of Walter's previous work on global continuation of steady water waves. The question can also be regarded as the compressible analog of the Maclaurin ellipsoids for incompressible rotating fluids. However, there are fundamental difficulties with the variational method mentioned above when it comes to proving continuation results. In particular, the non-convexity of the energy functional related to this problem makes it very difficult to prove uniqueness of minimizers (which may in fact be false in general). There is also no natural mechanism for continuous change of the minimizers when we continuously change the rotation speed.

Such a problem was partially resolved by our finding of Lichtenstein and Heilig's work. It is worth mentioning that Lichtenstein's original paper is in German, which is a language I cannot read. Walter, on the other hand, knows enough German to be able to confirm that Lichtenstein did have the basic result and idea for local perturbation. As we learned more about Lichtenstein and Heilig's work, it became clear to us that the lack of control on the total mass in their construction need to be remedied before we can globally continue the solution set to large rotation speed.

We did resolve this problem and upgraded Lichtenstein's theorem (Theorem \ref{thm: Lit}) to include a mass control (see \cite{strauss2017steady}):
\begin{theorem}\label{thm: AW1}
Let $\frac65<\gamma<2$, $\gamma\ne \frac43$, while other assumptions remain the same as in Theorem \ref{thm: Lit}. Then for each sufficiently small $\kappa$, there is a nonnegative compactly supported continuous function $\rho = \rho(\kappa)$ solving \eqref{eq: rot star} and $\int_{\real^3}\rho(x)~dx =\int_{\real^3}\rho_0(x)~dx$. The mapping $\kappa\mapsto \rho(\kappa)$ is continuous into a suitable function space, and $\rho(0)=\rho_0$.
\end{theorem}
Lichtenstein constructed his solutions by deforming the fluid domain and using an Ansatz for the rotating solutions. The main idea of our proof of Theorem \ref{thm: AW1} is to modify Lichtenstein's Ansatz so that a mass control will be enforced explicitly. We also need to make a technical change in the deformation map in order to help rigorously prove the estimates needed to apply the implicit function theorem in a suitable function space. Finally, the key new difficulty is in proving the linearized operator of the implicit function theorem has a trivial kernel. The modified construction respecting the mass control results in an integro-differential equation for functions in the kernel of the linearized operator, whereas Lichtenstein's construction only needs a vanishing theorem for an elliptic PDE. We found an interesting general condition for the kernel to be trivial (that even works for general equation of state different from a power law). To explain that condition, we define the function $M(a)$ to be the total physical mass of the radial solution to \eqref{eq: LE2} with center value $u(0)=a$. Our condition says the kernel is trivial if and only if $M'(a_0)\ne 0$, where $a_0 = u_0(0)$ is the center value corresponding to the Lane-Emden solution we perturb from. More informally, the condition means that the total mass of the non-rotating star has a genuine first order change as one changes the central density of the star. The curious omission of the case $\gamma=\frac43$ in Theorem \ref{thm: AW1} has to do with the fact that $M(a) = M(a_0)\left(a/a_0\right)^{\frac{3\gamma-4}{2\gamma-2}}$, which is a consequence of the scaling symmetry $u(x) \to \lambda^{\frac{2\gamma-2}{2-\gamma}}u(\lambda x)$ of \eqref{eq: LE2}. In particular, we see that $M'(a)=0$ when $\gamma=\frac43$. This is a pathological case, as all rescaled Lane-Emden solutions of different sizes have the same total mass.

In the same paper \cite{strauss2017steady}, we proved a similar theorem for the Vlasov-Poisson equation.
At about the same time, Jang and Makino \cite{jang2017slowly} studied local perturbations of the Lane-Emden equations without using an explicit Ansatz as Lichtenstein and we did. Their result does not contain a mass control, however. Jang met with Walter and me during the Spring 2017 semester program at ICERM (Brown University). The three of us decided to generalize the perturbative method to MHD-Euler-Poisson -- a model for rotating magnetic stars. In \cite{jang2019existence}, we proved the first existence result on rotating magnetic stars for small rotation and weak magnetic field. 

Walter and I thus participated in and witnessed a small revival of the perturbation methods for rotating stars. Walter's vision, however, has always been on the structure of large deviations from the non-rotating solution.

\section{Topological degree theory and global continuation}

There is a large established literature on global bifurcation and continuation method using topological degrees. As an example, we have the following global implicit function theorem.
\begin{theorem}   \label{thm: GIFT}
Let $X$ be a Banach space and let $U$ be an open subset of $X\times\real$.  
Let $F:U\to X$ be a $C^1$ mapping in the Fr\'echet sense.  
Let $(\xi_0, \kappa_0)\in U$ such that $F(\xi_0,\kappa_0)=0$.  
Assume that the linear operator $\frac{\pa F}{\pa\xi}(\xi_0,\kappa_0)$ is an isomorphism on $X$.  
Assume also that the mapping $(\xi.\kappa) \to F(\xi,\kappa)-\xi$ is compact from $U$ to $X$, and that $\frac{\pa F}{\pa\xi}(\xi,\kappa)-I\in L(X)$ is compact.  
Let $\S$ be the closure in $X\times\real$ of the solution set $\{(\xi,\kappa)\ \Big |\ F(\xi,\kappa)=0\}$. 
Let $\K$ be the connected component of $\S$ to which $(\xi_0,\kappa_0)$ belongs.   
Then one of the following three alternatives is valid. 
\begin{enumerate}[(i)]
\item $\K$ is unbounded in $X\times\real$.

\item $\K\backslash \{(\xi_0,\kappa_0)\}$ is connected.

\item $\K \cap \pa U \ne \emptyset$.  
\end{enumerate}
\end{theorem}

This is a standard theorem basically due to Rabinowitz. Theorem 3.2 in \cite{rabinowitz1971some} is
in the case that $U=X\times\real$ and under some extra structural assumption.  
A more general version also appears in Theorem II.6.1 of \cite{kielhofer2006bifurcation};   
its proof is easy to generalize to permit a general open set $U$.  
The case of a general open set $U$ also appears explicitly in \cite{alexander1976implicit}.  Roughly speaking, the suitable compactness assumptions allow one to define the Leray-Schauder degree for the mapping $F(\cdot, \kappa)$. If none of the alternatives holds, the component of the global solution set on either side of $(\xi_0,\kappa_0)$ will be compact. This will cause the Leray-Schauder degree to vanish for large $\kappa$. Homotopy invariance of the degree then implies that it will vanish for $\kappa$ close to $\kappa_0$. This contradicts the solution curve given by the usual local implicit function theorem, which would force the degree to be $\pm1$.

The three alternatives in the conclusion of Theorem \ref{thm: GIFT} are often termed the ``blow-up" case, the ``loop" case, and the ``meeting the boundary" case. At least one of these three cases must happen along the solution set. Such a conclusion shows that the solution set is ``large" in a certain sense, and is not just the local curve of the solutions given by the usual implicit function theorem.

We would like to apply the theory of global continuation via topological degree to the rotating star problem. Unfortunately, it is very difficult to establish the necessary compactness properties for the Lichtenstein type deformation constructions. Instead, we turned to another formulation of the problem. Our strategy is to invert the $\gamma-1$ power in \eqref{eq: rot star} to obtain a fixed point setup. In particular, we look for a function $\rho\in C_{loc}(\real^3)\cap L^1(\real^3)$
and a real number $\alpha$ such that
\eqn\label{eq: form 2}
\rho(x) = \left[\frac1{|\cdot|}*\rho(x)+\int_0^{r(x)}\omega^2(s)s~ds+\alpha\right]_+^{\frac1{\gamma-1}}
\eeqn
for all $x\in\real^3$. Here $f_+ = \max(f,0)$. This formulation has the advantage that the convolution with $\frac1{|x|}$ provides a simple source for compactness. Of course, every solution of \eqref{eq: form 2} is a solution to \eqref{eq: rot star}. However, \eqref{eq: form 2} misses many solutions of \eqref{eq: rot star}. The reason is that the physical problem doesn't require equality of the two sides of \eqref{eq: form 2} in the vacuum domain when $\rho(x)=0$. In particular, \eqref{eq: form 2} requires the terms in the square brackets to be non-positive when $\rho(x)=0$, whereas \eqref{eq: rot star} does not. This discrepancy becomes especially problematic when $\omega(r)$ does not decay at infinity, because in that case, the term involving $\omega(r)$ will be very large for large $r(x)$, causing the right hand side of \eqref{eq: form 2} to be outside of $L^1(\real^3)$. If we stay within the class of $\omega(r)$ with sufficient decay at $\infty$, however, \eqref{eq: form 2} is a viable approach. The set of $\omega(r)$ with rapid decay at infinity is already a large and interesting set of rotation profiles.  

To describe the result Walter and I obtained, let us define 
\eqn\label{def: F1}
\F_1(\rho,\alpha,\kappa) = \rho(\cdot)-\left[\frac1{|\cdot|}*\rho(\cdot)+\kappa^2 \int_0^{r(x)}\omega^2(s)s~ds+\alpha\right]_+^{\frac1{\gamma-1}}.
\eeqn
Here $\kappa$ is a parameter describing the intensity of rotation. Let $\rho_0$ be a radial Lane-Emden solution, and $\alpha_0$ be the number such that $\F_1(\rho_0,\alpha_0,0)=0$ (such a number is guaranteed to exist when $\rho_0$ is a Lane-Emden solution). Let $M = \int_{\real^3}\rho_0(x)~dx$, and define
\eqn
\F_2(\rho) = \int_{\real^3}\rho(x)~dx - M,
\eeqn
and the pair
\eqn
\F(\rho,\alpha,\kappa) = (\F_1(\rho,\alpha,\kappa),\F_2(\rho)).
\eeqn

We can now state our result as follows (see \cite{strauss2019rapidly} for details):
\begin{theorem}\label{thm: global cont}
Suppose $\frac65<\gamma<2$, $\gamma\ne \frac43$. Assume $\omega(r)$ has suitable decay as $r$ tends to infinity. There exists a set $\K$ of solutions to $\F(\rho,\alpha,\kappa)=0$ satsfying the following properties
\begin{enumerate}
\item $\K$ is a connected set in $C^1_c(\real^3)\times\real\times \real$.
\item $\K$ contains $(\rho_0,\alpha_0,0)$ together with a local curve of solutions around it.
\item If $\frac{4}{3}<\gamma<2$, then 
$$\sup \{|x| ~|~\rho(x)>0, (\rho,\alpha,\kappa)\in \K\}=\infty.$$
If $\frac65<\gamma<\frac43$, then 
either 
$$\sup \{|x| ~|~\rho(x)>0, (\rho,\alpha,\kappa)\in \K\}=\infty,$$
or
$$\sup \{\rho(x)~|~x\in\real^3, (\rho,\alpha,\kappa)\in \K \}= \infty.$$
\end{enumerate}
\end{theorem}
The last statement means that for the range of $\gamma$ we consider, either the supports become unbounded, or the densities become pointwise unbounded, as one continues along the solution set. Furthermore, if $\gamma>\frac43$, then the first alternative must hold. We thus construct, for the first time, a {\it connected set} of solutions that is {\it global}. Keeping the mass constant along the solution set turns out to be a key point of our methodology.

In the following, I list and compare several key features of the known existence results on rotating star solutions. In this table, the new result refers to Theorem \ref{thm: global cont}. The old results refer to previous theorems by other authors.
\begin{center}
\begin{tabular}{|p{1.2 in}|p{1.2 in }|p{1.2 in}|p{1.2 in}|}
\hline
& old results (variational)  & old results (perturbative) & new result (global continuation)\\
\hline
range of $\gamma$ & $(4/3,\infty)$ & $(6/5,2)$ &$(6/5,2)\setminus \{4/3\}$\\
\hline
mass constraint & yes & no & yes\\
\hline
allow large rotation & yes & no & yes\\
\hline
continuity of the solution set  & no & yes & yes\\
\hline
nature of singularity formulation &no & no &yes\\
\hline
\end{tabular}
\end{center}
\vspace{.5 in}

\section{Constructing Rapidly rotating stars}
In this section, I sketch the main ideas in the proof of Theorem \ref{thm: global cont}. We first put $\rho$ in the function space
$$ 
\C_s = \left\{  f:\real^3\to \real\ \Big|\ f \text{ is continuous, axisymmetric, even in }x_3, 
\text{ and } \|f\|_s <\infty\right\},  $$ 
where 
$$
\|f\|_s =: \sup_{x\in\real^3}\langle x\rangle^s|f(x)| <\infty.$$ 
The reason for this is simply to provide a straightforward definition of $\frac{1}{|x|}*\rho$ so that it decays properly at infinity. Let us now focus on the terms in the square brackets in \eqref{def: F1}. Assuming proper decay of $\omega(r)$, we get $\omega^2(r)r\in L^1(0,\infty)$. Denoting $j(x) = \int_0^{r(x)}\omega^2(s)s~ds$, and
$j_\infty = \lim_{r(x)\to\infty}j(x)$, we can rewrite the terms in the square brackets as
\eqn
\frac1{|\cdot|}*\rho(\cdot) + \kappa^2 (j(x)-j_\infty)+(\alpha+\kappa^2j_\infty)\to \alpha+\kappa^2j_\infty
\eeqn
as $r(x)\to \infty$. We see clearly that if $\alpha+\kappa^2j_\infty>0$, then any solution $\rho$ of $\F_1(\rho,\alpha,\kappa)=0$ will not be in $L^1(\real^3)$, because it tends to a positive constant as $r(x)\to \infty$. Thus the only way to set up this mapping consistently is to require $\alpha+\kappa^2j_\infty<0$. In fact, in this case we have $[\dots]<0$ for $x$ near infinity, thus any solution satisfying $\rho(x)=[\dots]_+^{1/(\gamma-1)}$ will be zero near infinity, thus is compactly supported. Nevertheless, to get quantitative estimates on the support and to close other estimates on the spaces, we actually need a little gap $\kappa^2j_\infty+\alpha<-\frac1N$. We can solve the global continuation problem with this $\frac1N$ gap, and finally patch up the solutions by letting $N\to\infty$. To highlight other ideas in the proof, let us ignore this technical gap and pretend the necessary estimates are available for $N=\infty$.

We now apply the global implicit function theorem, Theorem \ref{thm: GIFT}, by using $\xi = (\rho,\alpha)$, $X = \C_s\times \real$, $U = \{(\xi,\kappa) = (\rho,\alpha,\kappa)\in X\times\real ~|~ \kappa^2j_\infty+\alpha<0\}$. By the heuristic argument above, one can show that $\F$ maps $U$ into $X$ and is $C^1$. The needed compactness properties will following from the inverse Laplacian, or convolution with $\frac1{|x|}$. The next key condition to verify is that the linearized operator $\frac{\partial F}{\partial \xi}(\xi_0,0)$ is an isomorphism on $X$. This step is actually non-trivial, but the main difficulty was already resolved in our earlier paper \cite{strauss2017steady}. As is eluded to above, the key condition here turns out to depend only on properties of the non-rotating, radial, Lane-Emden solutions. The kernel is trivial if and only if $M'(a_0)\ne 0$, where $a_0 = u_0(0)$ is the center value corresponding to the Lane-Emden solution we perturb from, and $M(a)$ is the total mass of the solution to \eqref{eq: LE2} with central value $u(0)=a$. That the condition $M'(a_0)\ne 0$ is indeed satisfied for our range of $\gamma$ is then a simple consequence of the scaling symmetry of \eqref{eq: LE2}. 

The general Theorem \ref{thm: GIFT} then provides us with a solution set $\K$ of three alternatives labeled (i), (ii) and (iii). However, they are not specific enough to give us the results in Theorem \ref{thm: global cont}. First of all, we want to eliminate alternative (ii). This is an alternative that is often described as the ``loop" case. The possible existence of the loop case would significantly weaken our result, as it corresponds to a solution set with no blow-ups. 
To see that this cannot happen, we observe that a connected $\K\setminus \{(\rho_0,\alpha_0,0)\}$ must contain another non-rotating solution $(\rho_1,\alpha_1,0)\ne (\rho_0,\alpha_0,0)$. Study of this solution shows that it's a radial non-rotating Lane-Emden solution with different center density $\rho_1(0)\ne\rho_0(0)$, and the same total mass $\int_{\real^3}\rho_1(x)~dx=\int_{\real^3}\rho_0(x)~dx$. This contradicts the strict monotonicity of $M(a)$ since $M'(a)\ne 0$ for all $a$.  

We are left with alternatives (i) and (iii), known as the blow-up case, and the meeting the boundary case. Let us prove Theorem \ref{thm: global cont} by contradiction. Assume, therefore, for all $(\rho,\alpha,\kappa)\in \K$ that $\rho$ is uniformly bounded in $L^\infty(\real^3)$, and the support of $\rho$ is also uniformly bounded. We want to conclude from these assumptions that neither case (i) nor case (iii) in Theorem \ref{thm: GIFT} can happen, thus arriving at a contradiction. 

Suppose case (i) happens, then $\|\rho\|_{\C_s}+|\kappa|+|\alpha|$ is unbounded along $\mathcal K$. By our assumption of uniform bounds on $\rho$, however, $\|\rho\|_{\C_s}$ must be uniformly bounded. Thus $|\kappa|+|\alpha|$ is unbounded. Remember
\eqn\label{eq: arg 1}
\rho(x)  = \left[\frac1{|\cdot|}*\rho(x) + \kappa^2j(x) +\alpha \right]_+^{1/(\gamma-1)},
\eeqn
\eqn\label{eq: arg 2}
\int_{\real^3}\rho(x)~dx=M,
\eeqn
\eqn\label{eq: arg 3}
\kappa^2j_\infty+\alpha<0.
\eeqn
If $\kappa$ is bounded, then $\alpha\to-\infty$ by \eqref{eq: arg 3}. By the assumption on uniform $L^\infty$ bound on $\rho$ and uniform support bound, one can prove a uniform bound on the size of $\frac1{|\cdot|}*\rho(x)$. Thus by \eqref{eq: arg 1} $\rho\equiv 0$ as $\alpha\to-\infty$, violating the mass equation \eqref{eq: arg 2}. Thus $\kappa$ must be unbounded. Now if $\kappa\to\infty$, the terms in the square brackets in \eqref{eq: arg 1} will increase very rapidly due to the term $\kappa^2 j(x)$ (need to assume $j(x)$ is strictly increasing here, which amounts to suitable assumptions on $\omega(r)$), which will cause $\rho$ to be positive far outside, violating the common support on $\rho$.

The argument above shows case (i) cannot happen. Now assume case (iii) happens, so that $\kappa^2j_\infty+\alpha\to 0$ as one continues along the solution set. Since $\int_{\real^3}\rho ~dx = M$ and $\rho$ has a uniform support bound, we have a lower bound 
$$\frac{1}{|\cdot|}*\rho(x)\gtrsim \frac1{|x|}$$
as $r(x)\to\infty$. Thus
$$\left[\frac{1}{|\cdot|}*\rho + \kappa^2(j-j_\infty)+\kappa^2j_\infty+\alpha\right]>0$$
as $r(x)\to\infty$. To get the last inequality, we need $j(x)\to j_\infty$ sufficiently rapidly as $r(x)\to\infty$, which again amounts to suitable assumptions on $\omega(r)$. This implies $\rho=[\dots]_+^{1/(\gamma-1)}$ is positive when $r(x)$ is large, and again contradicts the uniform support bound on $\rho$.

The contradiction above shows that either the $L^\infty$ norm of $\rho$ blows up, or the support of $\rho$ blows up, as one moves along the solution set. To get the final refinement that the support of $\rho$ must blow up when $\frac43<\gamma<2$, one just need to get a uniform $L^\infty$ bound on $\rho$. We can start from the obvious $L^1$ bound on $\rho$, and use $L^p$ type estimates on $\frac1{|x|}*\rho$ and the equation $\rho(x)=\left[\frac1{|\cdot|}*\rho(x) + \kappa^2j(x) +\alpha \right]_+^{1/(\gamma-1)}$ to iteratively improve the the exponent $p$ until we eventually reach a uniform $L^\infty$ bound. This can be done when $\frac1{\gamma-1}$ is sufficiently low and the support of $\rho$ is uniformly bounded.

\section{Future directions}

Walter has always been spirited and explorative when it comes to extending the boundary of mathematical knowledge. The above discussion of global continuation of rotating stars is just one of the many examples where he takes a fresh new look on an age old problem, and offers wonderful novel insight into the structure of the problem.

There are many further questions one could ask with this new point of view on rotating star solutions. For instance, can one extend the range of $\gamma$ for these global continuation results to include $\gamma>2$, just like the solutions obtained by variational methods? Can one remove the decay assumption on $\omega(r)$ in the global continuation, and allow in particular, a constant rotation profile? Can one prove these results for the equation of state of white dwarf stars
$$p(\rho) = \int_0^{\rho^{1/3}}\frac{x^4}{\sqrt{1+x^2}}~dx$$
rather than just a power law? (We have already made significant progress on this problem.) Can one prove similar results for other models, such as the Vlasov-Poisson equation, or the general relativisitic Euler equations? What is the significance of these methods for numerical computations of rotating stars? Do these methods provide new insight on the problem of stability of rotating stars?

It is marvelous to see Walter continuing making his contribution on these interesting problems, and more questions to come.
\bibliographystyle{acm}
\bibliography{rotstarbiblio}

\end{document}